\documentclass[12pt]{amsart}

\usepackage{amsmath,amsthm,amsfonts,amscd,amssymb,euscript,enumerate}

\swapnumbers

\newtheorem{Lemma}{Lemma}
\newtheorem{Theorem}[Lemma]{Theorem}

\newtheorem*{Theorem*}{Lemma}
\theoremstyle{remark}
\newtheorem{Remark}[Lemma]{Remark}

\newtheorem*{Remark*}{Remark}
\theoremstyle{definition}
\newtheorem{Definition}[Lemma]{Definition}

\numberwithin{Lemma}{section}

\numberwithin{equation}{section}

\newcommand{\R}{\ensuremath{{\mathbb R}}}

\newcommand{\dc}{\ensuremath{C^M}}

\newcommand{\oka}{\ensuremath{\mathcal O}}

\newcommand{\stb}{{\EuScript E}}

\newcommand{\An}{{\EuScript A}}

\newcommand{\ideal}{{\mathcal I}}

\newcommand{\zeri}{{\mathcal Z}}

\newcommand{\gta}{{\mathfrak a}}

\begin{document}

\title[Global Denjoy-Carleman Nullstellensatz]{A global Nullstellensatz for ideals of Denjoy-Carleman functions}

\author[F. Acquistapace, F. Broglia]{Francesca Acquistapace, Fabrizio Broglia}

\address{Dipartimento di Matematica, Universit\`a degli Studi di Pisa, Largo Bruno Pontecorvo, 5, 56127 Pisa, Italy}

\email{acquistf@dm.unipi.it, broglia@dm.unipi.it}

\author[A. Nicoara]{Andreea Nicoara}

\address{Department of Mathematics, University of Pennsylvania, 209 South $33^{rd}$ St.,  Philadelphia, PA 19104, United States}

\email{anicoara@math.upenn.edu}

\thanks{The work presented in this paper was initiated while the authors were visiting the American Institute of Mathematics, Palo Alto, CA for a SQuaRE workshop. The authors are grateful to AIM for the hospitality and financial support.}



\begin{abstract}
We prove a Nullstellensatz result for global ideals of Denjoy-Carleman functions in both finitely generated and infinitely generated cases.
\end{abstract}

\maketitle


\section{Introduction}

Let $\dc(\R^n)$ be the ring of Denjoy-Carleman functions in some class $M$ defined on $\R^n$, and let $\gta$ be an ideal in this ring. The increasing sequence of positive real numbers $M=\{M_0, M_1,
M_2, \dots \}$ that determines the Denjoy-Carleman class is chosen so that the resulting ring $\dc(\R^n)$ is composed of quasi-analytic functions and so that the process of resolution of singularities applies to $\dc(\R^n);$ see  \cite{bm} and \cite{bmv}. The classical theory of the Denjoy-Carleman quasi-analytic functions may be found in \cite{rudin}, while \cite{thilliez} is a comprehensive recent survey. We will characterize here the ideal $\ideal(\zeri(\gta))$ of functions  vanishing on the variety of the ideal $\gta$  in terms  of the \L ojasiewicz radical, another notion of \L ojasiewicz radical computed on compact sets, and the saturation:

\begin{Theorem}\label{mainthm} Let $\gta$ be an ideal in $\dc(\R^n).$

\begin{enumerate}[(i)]

\item If $\gta$ is finitely generated, then $\ideal(\zeri(\gta)) =  \widetilde { {\sqrt[\text{\L}]{\gta}} },$ the saturation of the \L ojasiewicz radical of $\gta.$

\item $\ideal(\zeri(\gta)) = \left({\sqrt[\text{\L}]{\gta}}\right)_K,$ where $ \left({\sqrt[\text{\L}]{\gta}}\right)_K$ is the \L ojasiewicz radical of $\gta$ computed on all compact subsets of $\R^n.$

\end{enumerate}

\end{Theorem}

The Main Theorem will be derived from the following result proven in \cite{abn}: 

\begin{Theorem}\label{iberoamthm}

Let $\gta \subset \stb (\R^n)$ be a \L ojasiewicz ideal, where $\stb (\R^n)$ is the ring of smooth functions on $\R^n.$ Then the following hold:

\begin{enumerate}[(a)]

\item $   \widetilde{\sqrt[\text{\L}]{\gta}}=\overline{\sqrt[\text{\L}]{\gta}},$ i.e. the saturation of $\sqrt[\text{\L}]{\gta}$ equals its Whitney closure.


\item $\ideal(\zeri(\gta)) = \overline{\sqrt[\text{\L}]{\gta}}.$

\end{enumerate}

\end{Theorem}

\section{Definitions}

The statements of Theorems ~\ref{mainthm} and ~\ref{iberoamthm} contain the saturation, the Whitney closure, two notions of \L ojasiewicz radical as well as the notion of a \L ojasiewicz ideal. We define all of them now:

\begin{Definition}
The {\sl saturation} of  an ideal $\gta$ in $\dc(\R^n)$  is the ideal $$ \tilde \gta = \{ g\in \dc (\R^n) \, | \,\forall x \in \R^n \ g_x \in \gta \dc_x \},$$ where $\dc_x$ is the local ring of germs of Denjoy-Carleman functions in the class $M$ at the point $x.$ The {\sl saturation} of  an ideal $\gta$ in $ \stb (\R^n)$  is the ideal $$ \tilde \gta = \{ g\in \stb (\R^n) \, | \,\forall x \in \R^n \ g_x \in \gta \stb_x \},$$ where $\stb_x$ is the ring of germs of smooth functions in $n$ variables at $x.$
\end{Definition}

As shown in \cite{abn}, if $\gta$ is an ideal in $\stb (\R^n),$ then 
\begin{equation*}
\begin{split}
\tilde \gta &= \{g\in\stb (\R^n) \, | \,\forall \text{ compact }  K\subset \R^n \ \exists\,h\in \stb (\R^n)\:\text{s.t.}\:\zeri(h)\cap K=\emptyset \: \text{and} \:hg\in\gta\} \\&= \{g\in \stb (\R^n) \, | \,\forall\,x\in \R^n\,\exists\,h\in \stb (\R^n)\:\text{s.t.}\: h(x)\neq0 \: \text{and} \:hg\in\gta\},
\end{split}
 \end{equation*}
where $\zeri(h)$ is the zero set of $h.$

\begin{Definition}
Let the  algebra $ \stb  (\R^n)$ of  smooth
real-valued functions on $\R^n$ be endowed with the compact open topology. Let $\gta$ be an ideal in $ \stb  (\R^n).$ The {\sl Whitney closure} $\overline{\gta}$ is the closure of $\gta$ in this topology. By the Whitney Spectral Theorem, $\overline{\gta}=\gta^*,$ where $$ \gta^* = \{ g\in \stb (\R^n) \, | \,\forall x \in \R^n\ T_x g \in  T_x\gta \}.$$
\end{Definition}

\begin{Definition}\label{L radical}
The {\sl \L ojasiewicz radical}  of an ideal $\gta$ in $\stb(\R^n)$ is given
by  $$\sqrt[\text{\L}]{\gta}:=\{g\in\stb   (\R^n)\,  |  \,  \exists\,f\in\gta\:
\text{and} \: m\geq1\:\text{s.t.}\: f-g^{2m}\geq0 \,\, \forall x \in \R^n\}.$$ If $\gta$ is an ideal in $\dc(\R^n),$ then $$\sqrt[\text{\L}]{\gta}:=\{g\in\dc(\R^n)\,  |  \,  \exists\,f\in\gta\:
\text{and} \: m\geq1\:\text{s.t.}\: f-g^{2m}\geq0 \,\, \forall x \in \R^n\}.$$
\end{Definition}

\begin{Definition}\label{L radical compact}
The {\sl \L ojasiewicz radical computed on all compact subsets of $\R^n$}  of an ideal $\gta $ in $\dc(\R^n)$ is given
by  
\begin{equation*}
\begin{split}
 \left({\sqrt[\text{\L}]{\gta}}\right)_K :=\{ g\in\dc(\R^n) \, | \, \forall \text{ compact }  K\subset \R^n\, \exists\,f\in\gta\:\:
\text{and} \:\: m\geq1\:\text{s.t.}\: f-g^{2m}\geq0& \,\, \\ \forall x \in K&\}.
\end{split}
\end{equation*}
\end{Definition}

\begin{Definition}\label{Loj}
An ideal  $\gta \subset \stb (\R^n)$ is a {\sl  \L  ojasiewicz  ideal} if 

\begin{enumerate}
\item $\gta$ is generated by finitely many smooth functions $f_1,\ldots,f_l;$
\item $\gta$ contains an element $f$ with the property that for any compact $K\subset  M$ there exist a constant $c$ and an integer $m$ such that $|f(x)|\geq c\, d(x,  \zeri(\gta))^{m}$ on an open neighborhood of $K$.   
\end{enumerate}
\end{Definition}

\begin{Remark}
The element $f$ can be taken to be the sum of squares of the generators $f_1^2 +\dots +f_l^2.$
\end{Remark}

\section{Proof of Theorem~\ref{mainthm}}

\noindent {\bf Proof of part (i) of Theorem~\ref{mainthm}:} The ideal $\gta$ is finitely generated so $\gta = (f_1, \dots, f_l).$ In \cite{bm} Bierstone and Milman extended the resolution of singularities to the Denjoy-Carleman quasi-analytic classes. As a result, each of the generators $f_i$ of $\gta$ satisfies a \L ojasiewicz inequality; see Theorem 6.3 of \cite{bm}. Let $\hat \gta$ be the ideal in the ring $\stb (\R^n)$ generated by the elements of $\gta.$ Clearly, $\hat \gta$ is a \L ojasiewicz ideal since $f=f_1^2 +\dots +f_l^2$ satisfies part (2) of Definition~\ref{Loj}. Theorem~\ref{iberoamthm} then implies that $\ideal(\zeri(\hat\gta)) = \widetilde{\sqrt[\text{\L}]{\hat\gta}}.$ Since $\zeri(\gta) = \zeri (\hat \gta),$ it follows that $\ideal(\zeri(\gta)) = \ideal(\zeri(\hat\gta)) \cap \dc (\R^n) =  \widetilde{\sqrt[\text{\L}]{\hat\gta}} \cap \dc(\R^n).$ If we could show that $  {\sqrt[\text{\L}]{\gta}}  =\sqrt[\text{\L}]{\hat\gta} \cap \dc(\R^n),$ then it would follow
 that $ \widetilde { {\sqrt[\text{\L}]{\gta}} } = \widetilde{\sqrt[\text{\L}]{\hat\gta}} \cap \dc(\R^n),$ where saturation on the left-hand side is in $\dc(\R^n),$ while saturation on the right-hand side is in $\stb(\R^n).$ $  {\sqrt[\text{\L}]{\gta}}  =\sqrt[\text{\L}]{\hat\gta} \cap \dc(\R^n)$ is a consequence of the following lemma:

\begin{Lemma}\label{ra approx}
Let $\oka(\R^n)$ be the ring of real-analytic functions on $\R^n,$ and let $\An(\R^n)$ be any ring of functions on $\R^n$ such that $\oka(\R^n) \subset \An(\R^n) \subset \stb(\R^n).$ If $\gta$ is any ideal in $\An (\R^n)$ and if $\hat \gta$ is the ideal generated by the elements of $\gta$ in $\stb(\R^n),$ then $  {\sqrt[\text{\L}]{\gta}}  =\sqrt[\text{\L}]{\hat\gta} \cap \An(\R^n).$
\end{Lemma}

\begin{Remark}
It is not necessary for $\gta$ to be finitely generated.
\end{Remark}

\begin{proof}

Note that $  {\sqrt[\text{\L}]{\gta}}  \subset \sqrt[\text{\L}]{\hat\gta} \cap \An(\R^n)$ is trivially true, so the proof of the lemma reduces to showing the reverse inequality $\sqrt[\text{\L}]{\hat\gta} \cap \An(\R^n) \subset {\sqrt[\text{\L}]{\gta}} .$

\smallskip\noindent  Let $g  \in  \sqrt[\text{\L}]{\hat\gta} \cap  \An(\R^n),$
then there exist $f \in \hat\gta$ and  $m \geq 1$ such that $f \geq g^{2m}$ on
$\R^n.$  The  latter  implies  $f^2  \geq  g^{4m}$ on  $\R^n.$  Since  $f  \in
\hat\gta,$  there exist  $f_1,  \dots, f_k  \in  \gta,$ $h_1,  \dots, h_k  \in
\stb(\R^n)$ such that $f = h_1 f_1  + \cdots + h_k f_k.$ Therefore, $|f | \leq
|h_1| \ |f_1|  + \cdots + |h_k| \ |f_k|,  $ so $$ |f|^2 \leq  (|h_1| \ |f_1| +
\cdots + |h_k|  \ |f_k|)^2 \leq k \,  (|h_1|^2 \ |f_1|^2 + \cdots  + |h_k|^2 \
|f_k|^2).$$  If we  could bound  from  above the  coefficient functions  $h_1,
\dots, h_k$ by  functions in $\An(\R^n),$ we would be  done. Let $\epsilon >0$
be  given.  Since the  coefficient  functions  are  smooth, they  are  clearly
continuous, so there exist  global real-analytic functions $\tilde h_1, \dots,
\tilde h_k  \in \oka(\R^n)$ such that  $|\tilde h_i (x) -  h_i(x)| < \epsilon$
for all $ x  \in \R^n$ and all $i= 1, \dots, k.$  Therefore, $h_i (x) < \tilde
h_i(x) +  \epsilon$ for all  $ x \in \R^n$  and all $i= 1,  \dots, k.$ Note
that $\oka(\R^n) \subset \An(\R^n).$ We have thus shown that $$g^{4m} \leq f^2
<  k   \,  \left((\tilde  h_1^2+\epsilon)^2   \  f_1^2  +  \cdots   +  (\tilde
  h_k^2+\epsilon)^2  \  f_k^2\right),$$   where  $\tilde  f=k  \,\left((\tilde
  h_1^2+\epsilon)^2   \  f_1^2   +  \cdots   +  (\tilde   h_k^2+\epsilon)^2  \
  f_k^2\right) \in \gta$ as needed.

\end{proof}

We now let $\An(\R^n) = \dc (\R^n)$ noting that global real-analytic functions are contained in any Denjoy-Carleman ring of global functions of the type we are considering. We obtain $  {\sqrt[\text{\L}]{\gta}}  =\sqrt[\text{\L}]{\hat\gta} \cap \dc(\R^n),$ hence $ \widetilde { {\sqrt[\text{\L}]{\gta}} } = \widetilde{\sqrt[\text{\L}]{\hat\gta}} \cap \dc(\R^n).$\qed

\begin{Remark}
It should be noted here that part (i) of Theorem~\ref{mainthm} parallels the result obtained by Acquistapace, Broglia, and Fernando in \cite{abf} for an ideal of global real-analytic functions on a ${\mathcal C}$-analytic set.
\end{Remark}

\medskip
\noindent {\bf Proof of part  (ii) of Theorem~\ref{mainthm}:} The ideal $\gta$
here is not  necessarily finitely generated.  Note  that only
$\ideal(\zeri(\gta)) \subset \left({\sqrt[\text{\L}]{\gta}}\right)_K$ needs to
be     proven      as     $     \left({\sqrt[\text{\L}]{\gta}}\right)_K\subset
\ideal(\zeri(\gta))$  is   obvious.  Consider  $g   \in  \ideal(\zeri(\gta)).$
Let $K$
be  any compact subset  of $\R^n.$  We will  use Topological  Noetherianity, a
consequence  of  the  resolution  of  singularities  for  the  Denjoy-Carleman
classes, to  show there exist a  finite number of elements  of the original
ideal $f_1,  \dots, f_l  \in \gta$  and an open  set $U  \supset K$  such that
$  \zeri(\gta)\cap U =  \zeri\left(  (f_1, \dots,  f_l)
\right)\cap U.$ If $\gta = (0),$ there is nothing to be proven, so assume
$\gta \neq (0),$ and take  some
$f_1  \in  \gta$ not identically zero.  If there  exists  an  open  set  $U  \supset K$  such  that
$\zeri((f_1))= \zeri (\gta)$ on $U,$ we are done with $l =1;$ otherwise, there
is some $f_2 \in \gta$  such that $\zeri((f_1, f_2)) \subsetneq \zeri((f_1)).$
Inductively,  we  will  have chosen  $f_1,  \dots,  f_k  \in \gta$  such  that
$\zeri((f_1, f_2, \dots, f_k))  \subsetneq \zeri((f_1, f_2, \dots, f_{k-1})).$
Hence we have a sequence of finitely generated ideals $(f_1) \subset (f_1,f_2)
\subset \ldots \subset (f_1, \ldots, f_k) \subset \dots$.

Topological Noetherianity, Theorem 6.1  in \cite{bm} (see also Theorem  3.1   in  \cite{bmv}), guarantees that the
sequence of corresponding zero sets stabilizes, i.e. there exists some $k$ and some open set $U \supset K$ such
that $\zeri((f_1, f_2, \dots, f_k)) = \zeri(\gta)$ on $U.$ Set $l=k.$  We   have  obtained  $
\zeri(\gta)\cap U =  \zeri\left(  (f_1, \dots, f_l)  \right)\cap U.$ We now
recall the following lemma from \cite{abn}:

\begin{Lemma}\label{L disug}
Let  $\gta$ be  a  \L  ojasiewicz ideal  generated  by $f_1,\ldots,f_l$  and
$f=f_1^2+\dots +f_l^2$.   Let $g\in \stb (\R^n)$ be  such that $\zeri(g)\supset
\zeri(f) = \zeri(\gta)$. Then for  any compact set $K\subset \R^n,$ there exist
a constant  $c$ and a positive  integer $m$  such that $g^{2m}  \leq cf$ on  an open
neighborhood  of $K$.  In  particular, there  exist  an integer  $m$ and  an
element $a\in \gta$  such that $g^{2m} \leq |a|$ on  an open neighborhood of
$K$.
\end{Lemma}

\medskip Its proof from \cite{abn} carries over if instead of $\R^n$ we restrict to an open set $U \subset \R^n$ provided that all compact sets considered are subsets of $U.$ The ideal generated by $f_1, \dots, f_l$ in $\stb (\R^n)$ is \L ojasiewicz and since $g \in \ideal(\zeri(\gta)),$ it follows $\zeri(g)\Big|_U \supset \zeri\left(  (f_1, \dots, f_l)  \right)\Big|_U.$ By the lemma, there exist a constant  $c$ and a positive  integer $m$  such that $g^{2m}  \leq cf$ on some set $\tilde U,$ where $K \subset \tilde U \subset U.$ Therefore, $g \in  \left({\sqrt[\text{\L}]{\gta}}\right)_K .$ \qed

\begin{Remark}
Unlike  $\sqrt[\text{\L}]{\cdot},$ $ \left({\sqrt[\text{\L}]{\cdot}}\right)_K$ already contains in its definition a saturation operation. In other words, for any ideal $\gta$ in $\dc(\R^n),$ $ \left({\sqrt[\text{\L}]{\gta}}\right)_K$ is already saturated, while $\sqrt[\text{\L}]{\gta}$ might not be.
\end{Remark}

\bibliographystyle{plain}
\bibliography{GlobalNSS}

\end{document}